\documentclass[11pt]{amsart}

\usepackage{amsmath,amssymb}

\begin{document}

\title[a priori]{  A priori bound on the velocity in  \\
 axially symmetric  Navier-Stokes equations}
\thanks{}
\date{Sept. 2013; revised Feb. 2015}

\vspace{ -1\baselineskip}

\author{ Zhen Lei, Esteban A Navas and Qi S. Zhang}

\address{
Z. L.: School of Mathematics, Fudan University, Shanghai 200433, P.
R. China.}

\address{E. A. N. and Q.S.Z.: Department of Mathematics,  University of California,
Riverside, CA 92521}

\numberwithin{equation}{section}
\allowdisplaybreaks

\newtheorem{theorem}{Theorem}
\newtheorem{proposition}{Proposition}
\newtheorem{lemma}{Lemma}
\newtheorem{corollary}{Corollary}
\newtheorem{remark}{Remark}[section]
\newtheorem{acknowledgement}{Acknowledgement}
\numberwithin{theorem}{section}
\numberwithin{lemma}{section}
\numberwithin{proposition}{section}
\numberwithin{equation}{section}
\def\R{\bf R}
\def\eps{\varepsilon}
\def\b{\bar}
\def\al{\aligned}
\def\eal{\endaligned}
\def\be{\begin{equation}}
\def\ee{\end{equation}}
\def\lab{\label}
\def\nn{\nonumber}
\def\~{\widetilde}
\def\->{\overrightarrow}
\def\bar{\overline}

\newcommand{\reals}{\mathbb{R}}
\newcommand{\Div}{\textrm{div }}
\newcommand{\Curl}{\textrm{curl }}
\newcommand{\supp}{\textrm{supp }}
\newcommand{\V}{\widetilde{v}}
\newcommand{\Omegabar}{\overline{\Omega}}


\newcommand{\Z}{{\mathbb Z}}
\newcommand{\N}{{\mathbb N}}
\newcommand{\F}{{\mathbb F}}
\newcommand{\C}{{\mathbb C}}
\newcommand{\Q}{{\mathbb Q}}
\newcommand{\T}{{\mathbb T}}
\newcommand{\barp}{{\overline{P}}}
\newcommand{\barpt}{{\overline{P}_{\frac{1}{2},\frac{16}{3},1}}}
\newcommand{\barpf}{{\overline{P}_{\frac{1}{4},\frac{64}{9},1}}}
\newcommand{\bary}{{\overline{y}}}
\newcommand{\barg}{{\overline{\nabla}}}
\newcommand{\barc}{{\overline{C}}}
\newcommand{\bark}{{\overline{K}}}
\newcommand{\baro}{{\overline{\Omega}}}
\newcommand{\kbp}{{\overline{K}^{p'}(b, \overline{C}_{1,4,1})}}
\newcommand{\kbpf}{{\overline{K}^{p'}(b, \overline{C}_{\frac{1}{4},\frac{64}{9}, 1})}}
\newcommand{\vgn}{{||V||^{\frac{1}{1-\lambda}}_{L^{(2^+)^+}(\overline{P}_{1,4,1})}}}
\newcommand{\vp}{{\frac{1}{1-\lambda}}}
\newcommand{\vlp}{{(2^+)^+}}
\newcommand{\rt}{{\tilde{r}}}
\newcommand{\ttheta}{{\tilde{\theta}}}
\newcommand{\z}{{\tilde{z}}}
\newcommand{\x}{{\tilde{x}}}
\newcommand{\vb}{{\overline{v}}}
\newcommand{\ball}{{B_{\frac{r_0}{4}}}}
\newcommand{\ballt}{{B_{2r_0}}}
\newcommand{\omt}{{\frac{C}{r_0^{3.5^+}}}}

 \smallskip

\begin{abstract}
Let $v$ be the velocity of Leray-Hopf solutions to the axially
symmetric three-dimensional Navier-Stokes equations. Under suitable conditions for
initial values,
we prove the following a priori bound
\[
|v(x, t)| \le \frac{C |\ln r|^{1/2}}{r^2},  \qquad 0< r \le 1/2,
\]where  $r$ is the distance
from $x$ to the $z$ axis, and
$C$ is a constant depending only on the initial value.

This provides a pointwise upper bound (worst case scenario) for
possible singularities while the recent papers \cite{CSTY2} and
\cite{KNSS} gave a lower bound. The gap is polynomial order $1$
modulo a half log term.

\end{abstract}
\maketitle

\maketitle
\tableofcontents

\section{Introduction}

In this paper, we prove, under suitable initial condition, that the flow speed in the axially
symmetric incompressible flow has an a priori bound
which is proportional to  the inverse square of the distance to the
rotational axis, modulo a half logarithmic term.  In order to present the result precisely, let us first recall the basic
set ups.
In
Cartesian coordinates,
the incompressible Navier-Stokes equations are
\begin{align*}
\Delta v-(v\cdot\nabla)v-\nabla p-\partial_t v=0,\ \Div v=0,
\end{align*}
where
$v=(v_1(x,t),v_2(x,t),v_3(x,t)):\reals^3\times[0,T]\rightarrow\reals^3$
is the velocity field and
$p=p(x,t):\reals^3\times[0,T]\rightarrow\reals$ is the pressure.
In cylindrical coordinates $r,\theta, z$ with
$(x_1,x_2,x_3)=(r\cos\theta,r\sin\theta,z)$,  axially symmetric
solutions are of the form
\begin{align*}
v(x,t)=v_r(r,z,t)\->{e_r}+v_{\theta}(r,z,t)
\->{e_{\theta}}+v_z(r,z,t)\->{e_z}.
\end{align*}
The components $v_r,v_{\theta},v_z$
 are all independent of the angle of rotation $\theta$.
 Here $\->{e_r},\->{e_{\theta}},\->{e_z}$ are the basis vectors for
 $\reals^3$ given by
\begin{align*}
\->{e_r}= \Big (\frac{x_1}{r},\frac{x_2}{r},0 \Big ),
\ \->{e_{\theta}}= \Big (\frac{-x_2}{r},\frac{x_1}{r},0 \Big ),
\ \->{e_z}=(0,0,1).
\end{align*}

It is well known   (see
\cite{CSTY1} for example) that $v_r$, $v_z$ and $v_\theta$ satisfy the equations
\begin{align}
\lab{eqasns}
\begin{cases}
   \big (\Delta-\frac{1}{r^2} \big )
v_r-(b\cdot\nabla)v_r+\frac{v_{\theta}^2}{r}-\frac{\partial
p}{\partial r}-\frac{\partial v_r}{\partial t}=0,\\
   \big   (\Delta-\frac{1}{r^2}  \big
)v_{\theta}-(b\cdot\nabla)v_{\theta}-\frac{v_{\theta}v_r}{r}-\frac
{\partial v_{\theta}}{\partial t}=0,\\
 \Delta v_z-(b\cdot\nabla)v_z-\frac{\partial p} {\partial
z}-\frac{\partial v_z}{\partial t}=0,\\
 \frac{1}{r}\frac{\partial (rv_r)}{\partial r} +\frac{\partial
v_z}{\partial z}=0,
\end{cases}
\end{align}
 where $b(x,t) =(v_r,0,v_z)$ and
the last equation is the divergence-free
 condition. Here, $\Delta$ is the cylindrical  scalar Laplacian
and $\nabla$ is the cylindrical gradient field:
 \begin{align*}
 \Delta =\frac{\partial^2}{\partial r^2}+
\frac{1}{r}\frac{\partial}{\partial r}+\frac{1}{r^2}\frac{\partial
^2}{\partial \theta ^2}+\frac{\partial^2}{\partial
 z^2},\ \
 \nabla =  \Big  (\frac{\partial }{\partial r},\frac{1}{r}\frac{\partial
}{\partial \theta},\frac{\partial }{\partial  z}  \Big  ).
\end{align*}
Observe that the equation for $v_{\theta}$ does not explicitly depend on the
pressure. Let $\Gamma=rv_{\theta}$, one sees that the function
$\Gamma$ satisfies
\begin{equation}\label{Gamma/vtheta}
\Delta \Gamma -(b\cdot\nabla)\Gamma-\frac{2}{r}\frac{\partial
\Gamma}{\partial r}-\frac{\partial\Gamma}{\partial t}=0,\ \Div b=0.
\end{equation}
Recall that the vorticity  $\omega=\Curl v$ for axially
symmetric solutions
 \begin{align*}
 \omega(x,t)=\omega_r\->{e_r}+\omega_{\theta}
\->{e_{\theta}}+\omega_z\->{e_z}
\end{align*} is given by
\begin{equation}
\label{curlformulas} \omega_r =-\frac{\partial v_{\theta}}{\partial
z},\ \omega_ {\theta} =\frac{\partial v_r}{\partial
z}-\frac{\partial v_z}{\partial r},\ \omega_z=\frac{\partial v_
{\theta}}{\partial r}+\frac{v_{\theta}}{r}.
\end{equation}
The equations of vorticity $\omega= \Curl v$ in
cylindrical form are  (again, see \cite{CSTY1} for example):
 \begin{align*}
\begin{cases}
  \big (\Delta-\frac{1}{r^2} \big
)\omega_r-(b\cdot\nabla)\omega_r+\omega_r\frac {\partial
v_r}{\partial r}+\omega_z\frac{\partial v_r}{\partial
z}-\frac{\partial
\omega_r}{\partial t}=0,\\
   \big  (\Delta-\frac{1}{r^2}  \big
)\omega_{\theta}-(b\cdot\nabla)\omega_{\theta}+2\frac{v_{\theta}}
{r}\frac{\partial v_{\theta}}{\partial
z}+\omega_{\theta}\frac{v_r}{r}-\frac{\partial
\omega_{\theta}}{\partial t}=0,\\
 \Delta\omega_{z}-(b\cdot\nabla)\omega_{z}+\omega_z\frac{\partial
v_z}{\partial z}+\omega_{r}\frac{\partial v_z}{\partial
r}-\frac{\partial \omega_{z}}{\partial t}=0.
\end{cases}
\end{align*}
 Define
$\displaystyle\Omega=\frac{\omega_{\theta}}{r}$, then we have that
$\Omega$ satisfies
\begin{equation}\label{Omega/omegatheta}
\Delta \Omega -(b\cdot\nabla)\Omega+\frac{2}{r}\frac{\partial
\Omega}{\partial r}-\frac{\partial\Omega}{\partial
t}+\frac{2v_{\theta}}{r^2}\frac{\partial v_{\theta}}{\partial z}=0,\
\Div b=0. \end{equation}

 If the swirl $v_\theta=0$, then it is
known for long time (see O. A. Ladyzhenskaya
\cite{L}, M. R. Uchoviskii and B. I. Yudovich \cite{UY}), that finite energy solutions
to (\ref{eqasns}) are smooth for all time. See also the paper by S.
Leonardi, J. Malek, J. Necas,  and  M. Pokorny \cite{LMNP}).

In the presence of swirl, it is not known in general if finite energy solutions
blow up in finite time.   However a lower bound for the possible blow up rate is known by the recent results of
C.-C. Chen, R. M. Strain, T.-P.Tsai,  and  H.-T. Yau in \cite{CSTY1},
\cite{CSTY2},  G. Koch, N. Nadirashvili, G. Seregin,  and  V.
Sverak in \cite{KNSS}. See also the work by G. Seregin  and
V. Sverak \cite{SS} for a localized version. These authors prove that if
\be
\lab{v<1/r}
|v(x, t)| \le
\frac{C}{r},
\ee then solutions are smooth for all time. Here $C$ is any positive constant.  Solutions satisfying this bound are
often refereed to as type I solutions. One reason for this name is that the bound scales the same way as the reciprocal of the distance.
Their result can be rephrased as: type I solutions are regular.
See also the paper by two of us  \cite{LZ11}, \cite{LZ11-2} on further results in this direction.
Without knowing if blow up happens in general, it is desirable to find an upper
bound for the growth of velocity. It is expected that the solutions are
smooth away from the axis, with  certain  growing bound
when approaching the axis. Our Theorem \ref{ourtheoremV} confirms this
intuitive idea. Although it did not give the bound
(\ref{v<1/r}) which is required for smoothness, it
reveals the exact gap between what we have and what we need.

 This seems to be the first pointwise bound for the speed
(velocity) for the axially symmetric Navier-Stokes equation.
  We mention that
a less accurate a priori upper bound for the vorticity has been found in
\cite{BZ:1}.

The following are some notations to be frequently used. We use $x=(x_1,x_2,x_3)$ to denote
a point in $\reals^3$ for rectangular coordinates, and in the
cylindrical system we use $r=\sqrt{x_1^2+x_2^2}$,
$\theta=\tan^{-1}\frac{x_2}{x_1}$, $z=x_3$. Also $x'=(x_1, x_2, 0)$.

Here is the main result of the paper.

\begin{theorem} \label{ourtheoremV}
 Suppose  $v$ is a smooth,
axially symmetric solution of the three-dimensional Navier-Stokes
equations in $\reals^3\times(-T,0)$ with initial data
$v_0=v(\cdot,-T)\in L^2(\reals^3)$. Assume further
$rv_{0,{\theta}}\in L^{\infty}(\reals^3)$ and let $R= \min\{1,
\sqrt{T/2} \}$.

 Then for all $(x, t) \in \reals^3 \times (-R^2, 0)$, it holds

\[
|v_r(x, t)| + |v_z(x, t)| \le \frac{C \sqrt{|\ln r|}}{r^2},  \qquad 0< r
\le \min \{1/2, R \}.
\]Here  $r$ is the distance
from $x$ to the $z$ axis, and $C$ is a constant depending only on the initial data.

\end{theorem}

The proof of the theorem is based on the following pointwise bound on the vorticity.

\begin{theorem} \label{ourtheorem}
 Suppose  $v$ is a smooth,
axially symmetric solution of the three-dimensional Navier-Stokes
equations in $\reals^3\times(-T,0)$ with initial data
$v_0=v(\cdot,-T)\in L^2(\reals^3)$, and $\omega$ is the vorticity.
Assume further, $rv_{0,{\theta}}\in L^{\infty}(\reals^3)$ and let
$R= \min\{1, \sqrt{T/2} \}$. Then the following a priori estimate
holds.

 There is a constant $C$,
depending only on the initial data, such that the following holds for all  $(x, t)
\in \reals^3 \times (-R^2, 0)$ with $r=|x'| \in (0, \min\{ R,
\frac{1}{2} \})$:
\[
\al |\omega_{\theta}(x, t)| &\leq\frac{C \ln (1/r)}{r^{7/2}}     \left[
\sup_{s \in [t-r^2, t]} \left(  \int_{B(x, 4 r)} (v^2_r + v^2_z)(y,
s) dy \right)^{1/2} + r^{1/2} (\|r v_{0, \theta}
\|_{L^{\infty}(\reals^3)}+1)
   \right]^{2}  \\
&\hspace{3cm}  \times   \left[ \left(\int^t_{t-r^2} \int_{B(x, 4
r)} \omega^2_{\theta} (y, s) dyds \right)^{1/2}+ {r}^{1/2} (\|r
v_{0, \theta} \|_{L^{\infty}(\reals^3)}+1) \right]. \eal
\]
\end{theorem}

\begin{remark}
\rm We assume smoothness of the solution only for technical
simplicity. One can use standard approximation methods to treat
the weak solution case. In fact, if $v$ is a (suitable) Leray-Hopf
solution, then it is smooth except possibly on the $z$ axis (c.f.
\cite{CKN:1}). Also, the bound on $\omega_\theta$ is scaling
invariant. Similar bounds can also be proven for the other two
components of the vorticity $\omega_r$ and $\omega_z$. But we will
not do this here.
\end{remark}

Here we mention a number of related papers  on axially symmetric Navier-Stokes equations.  J. Neustupa  and
M. Pokorny \cite{NP} proved that the regularity of one component
(either $v_r$ or $v_{\theta}$) implies regularity of the other
components of the solution. Also proving regularity is the work of
Q. Jiu  and  Z. Xin \cite{JX} under an assumption of sufficiently small
zero-dimension scaled norms.  D. Chae  and  J. Lee \cite{CL}  also proved
regularity results assuming finiteness of another certain
zero-dimensional integral.  G. Tian  and  Z.
Xin \cite{TX}  constructed a family of singular axially
symmetric solutions with singular initial data.
T. Hou  and  C. Li \cite{HL} found a special class of global smooth
solutions. See also a recent extension: T. Hou, Z. Lei  and  C. Li
\cite{HLL}.

Let us outline the proof of the theorem. The starting point is the
a priori bound for the rotational component of the velocity: $r |v_\theta (\cdot, t)|  \in L^\infty$. A proof of this fact
can be found in \cite{CL} Section 3, Proposition 1, for example.
The first ingredient is the observation
that the basic energy estimate is critical when
localized in a dyadic ball which is away from the symmetric axis. This
enables us to perform a kind of dimension reduction argument
and apply two-dimensional Sobolev imbedding inequalities. The second
ingredient is a new estimate on the oscillation of the angular stream function
in a dyadic ball which is based on dimension reduction and the structure of the equations.
Then we apply the standard Moser iteration method
and two-dimensional Sobolev imbedding inequalities to get an upper
bound for the quantity $\Omega=w_\theta/r$, based on the evolution equation of $\Omega$ in
\eqref{Omega/omegatheta}. The
third ingredient is a novel use of the localized Biot-Savart law.
We use axis-symmetry to show that $L^2$ integrals of velocity in
small dyadic regions are smaller than usual. This fact and the a
priori bound on $\omega^\theta$ implies the point-wise bound on
$|v^r|+|v^z|$.

Throughout this paper we use $C$ to denote an absolute positive
constant. When $C$ depends on $p$, we use the notation $C_p$. The
meanings of $C$ and $C_p$ may change from line to line.

The remainder of the paper is organized as follows.
 In  Section 2,  we prove Theorem \ref{ourtheorem} the a priori bound on
$\omega_\theta$.  In Section 3 we prove
Theorem \ref{ourtheoremV}.  Some additional results are given in the appendix.

\section{A priori bound for $\omega_{\theta}$}

In this section, we will prove Theorem \ref{ourtheorem}.

Let $(x, t)$ be the point in
the statement of the theorem. For simplicity
we take $t=0$ and $x_3=0$.  During the proof, it is convenient to replace the three dimensional ball $B(x, 4 r)$ by comparable cylindrical
type regions. The reason is that a cylindrical region has a fixed profile in
the $e_r, e_z$ plane. This feature allows us to reduce much computations to 2 dimensional setting.

 So let us introduce a few more notations. Let $R>0$, $S>0$, and $0<A<B$ be
constants.

Denote
\be
\lab{R3tong}
C_{A R, B R}=\{(x_1,x_2,x_3) :  AR\leq r\leq BR,\ 0\leq\theta\leq
2\pi,\ |x_3|\leq BR\}\subset\reals^3
\ee
 to be the hollowed out cylinder
centered at the origin, with inner radius $AR$, outer radius $BR$,
and height extending up and down $BR$ units for a total height of
$2BR$. If $R=1$, we will write $C_{A, B}$ in place of $C_{A 1, B 1}$.

 Denote $P_{AR,BR,SR}$ to be the parabolic region
\be
\lab{paratong}
P_{AR,BR,SR}=C_{A R, B R}\times (-S^2R^2,0).
\ee If $R=1$, we will use
$P_{A, B, S}$ to denote $P_{A1,B1,S1}$.

Proving the theorem is equivalent to finding
priori bound for $\omega_\theta$ in the region
$P_{\frac{k}{2}, 2k,\frac{3k}{4}}$ with $0<k<\min\{1, \sqrt{T/2} \}$. We will use
the scaling property of the Navier Stokes equation to shift the
consideration to the cube $P_{\frac{1}{2}, 2,\frac{3}{4}}$.
 We recall that scaling of
the equations now; the pair $(v(x,t),p(x,t))$ is a solution to the
system, if and only if for any $k>0$ the re-scaled pair $( \V
(x,t),\widetilde{p}(x,t) )$ is also a solution, where $
 \V(x,t)=kv(kx,k^2t),$ $\widetilde{p}(x,t)=k^2p(kx,k^2t).
$ Thus, if $(v,p)$ is a solution to the axially symmetric
Navier-Stokes equations for $(x,t)\in P_{k,4k,k}$, then
$(\V(\widetilde{x},\widetilde{t}),
\widetilde{p}(\widetilde{x},\widetilde{t}))$ is a solution to the
equation in the variables
$\widetilde{x}=\frac{x}{k},\widetilde{t}=\frac{t}{k^2}$ when
$(\widetilde{x},\widetilde{t})\in P_{1,4,1}$. We note here how
certain quantities scale or change due to the above. Here, D is any
domain in $\reals^3$ and $kD=\{x:x=ky,\ y\in D\}$:
\begin{align*}
&r=\sqrt{x_1^2+x_2^2}:\quad \widetilde{r}= \sqrt{   \left
(\frac{x_1}{k}    \right  )^2+   \left (\frac{x_2}{k}   \right  )^2}
=\frac{r}{k}
\\
\|v(x,&t)\|_{L^2(kD\times(-(kR)^2,0))}:\\
& \|\~{v}(\~{x},\~{t})\|_{L^2(D\times(-R^2,0))} =  \Big
(\int_{-R^2}^0\int_D|\~{v}(\~{x},\~ {t})|^2d\~{x}d\~{t}  \Big
)^{\frac{1}{2}} \\
&=\Big (\int_{-(kR)^2}^0\int_{kD}|kv(x,t)|^2\frac{1}{k^5}dxdt  \Big
)^{\frac{1}{2}} =
\frac{1}{k^{\frac{3}{2}}}\|v(x,t)\|_{L^2(kD\times(-(kR)^2,0))}.
\end{align*}
\begin{align*}
b(x,&t)= (v_r,0,v_z): \\
&\widetilde{b}(x,t)=(kv_r(kx,k^2t),0,kv_z(kx,k^2t)) =kb(kx,k^2t),\
(x,t)\in P_{k,4k,k}\\  &\Rightarrow
\widetilde{b}(\widetilde{x},\widetilde{t})=kb(x,t).
\end{align*}
\begin{align*}
\|b(x,&t)\|_{L^{\infty}(-(kR)^2,0;L^2(kD))}: \\
&\|\~{b}(\~{x},\~{t})\|_{L^{\infty}(-R^2,0;L^2(D))} =\sup_{-R^2\leq
\~{t}<0}  \Big (\int_{D}|\~{b}(\~{x},\~ {t})|^2d\~{x} \Big
)^{\frac{1}{2}} \\
& =\sup_{-(kR)^2\leq t<0}   \Big
(\int_{kD}|kb(x,t)|^2\frac{1}{k^3}dx   \Big  )^{\frac{1}{2}} =
\frac{1}{k^{\frac{1}{2}}}\|b(x,t)
\|_{L^{\infty}{(-(kR)^2,0;L^2(kD))}}.
\end{align*}
\begin{align*}
\omega(x,&t): \quad\widetilde{\omega}(x,t)= k^2\omega(kx,k^2t),\
(x,t)\in P_{1,4,1} \Rightarrow
 \widetilde{\omega}(\widetilde{x},\widetilde{t})=k^2\omega(x,t)\\
\|\omega(x,&t)\|_{L^2(kD\times(-(kR)^2,0))}: \\ \!
 &|\~{\omega}(\~{x},\~{t})\|_{L^2(D\times \!
(-R^2,0))} \! =  \! \Big (\int_{-R^2}^0  \! \int_D \!
|\~{\omega}(\~{x},\~{t})|^2d\~{x}d\~{t}  \Big )^{\frac{1}{2}}
\\
& =     \Big
(\int_{-(kR)^2}^0\int_{kD}|k^2\omega(x,t)|^2\frac{1}{k^5}dxdt \Big
)^{\frac{1}{2}} =
\frac{1}{k^{\frac{1}{2}}}\|\omega(x,t)\|_{L^2(kD\times(-(kR)^2,0))}.
\end{align*}
One can also show that
$\widetilde{\Gamma}(\~{x},\~{t})=\widetilde{r}
\widetilde{v_{\theta}}(\~{x},\~{t})$ is a solution to
(\ref{Gamma/vtheta}) and $\widetilde{\Omega}(\~{x},\~{t})
=\frac{\widetilde{\omega_{\theta}}(\~{x},\~{t})}{\widetilde{r}}$ is
a solution to (\ref{Omega/omegatheta}) in the variables
$(\~{x},\~{t})\in P_{1,4,1}$. We will do most of our computations on
scaled cylinders.

Since $r v_\theta$ is scaling invariant, using the following
result, we know that $\~ r \~ v_\theta$ is uniformly bounded for
all time.

\begin{proposition} (\cite{CL}) and \cite{NP})
\label{chaeandleeprop}Suppose  $v$ is a smooth, axially symmetric
solution of the three-dimensional Navier-Stokes equations with
initial data $v_0\in L^2(\reals^3)$.  If $rv_{0,{\theta}}\in
L^p(\reals^3)$, then $rv_{\theta}\in L^{\infty}(0,T;L^p(\reals^3))$.
In particular, if $p=\infty$,
\begin{align*}
|v_{\theta}(x,t)|\leq \frac{\| rv_{0, \theta}\|
_{L^{\infty}(\reals^3)}}{\sqrt{x_1^2+x_2^2}}.
\end{align*}
\end{proposition}

 \smallskip

\noindent \textbf{Proof of Theorem \ref{ourtheorem}.}

During the proof, we are going to drop the ``tilde" notation for all
relevant quantities over a time when computations take place on the
scaled cylinders.  By the end, we will scale down to the original
solution. Although this scaling seems merely a technical move that
simplifies the computation, it is actually a key step that allows us to
care out a dimension reduction argument mentioned earlier. In the
region $P_{1,4,1}$ we do our analysis on (\ref{Omega/omegatheta}):
  \begin{align*}
\Delta \Omega -(b\cdot\nabla)\Omega+\frac{2}{r}\frac{\partial
\Omega}{\partial r}-\frac{\partial\Omega}{\partial
t}+\frac{2v_{\theta}}{r^2}\frac{\partial v_{\theta}}{\partial
z}=0,\ \Div b=0.
\end{align*}
A flow chart for the argument to prove  Theorem
\ref{ourtheorem} is as follows:

Step 1:  Energy Estimates by a refined cut-off function.

Step 2: Estimate drift term $(b\cdot\nabla)\Omega$  using methods
similar to \cite{Z}. Use dimension reduction. Note this term is more
singular than that allowed by standard theory.

Step 3: Estimate a term involving the cut-off.

Step 4: Estimate the term involving the directional
derivative $\partial_r$ using a method similar to that in \cite{CSTY1}.

Step 5: Estimate the inhomogeneous term utilizing the bound in
Proposition \ref {chaeandleeprop}.

Step 6: $L^2-L^{\infty}$ Estimate on Solutions to
(\ref{Omega/omegatheta}) via Moser's Iteration. Use dimension
reduction.

$L^2-L^{\infty}$ Estimate on $\omega_{\theta}$ via re-scaling.

\smallskip

\noindent
\textbf{Energy Estimates:}
\\
 \textbf{Step 1:} We use a revised cut-off function
and the equation to obtain inequality (\ref{T1T2T3T4}) below.

 Note that
\begin{equation}
\label{Lambda}
\Lambda \equiv
\| v_{\theta}\| _{L^{\infty}(P_{1,4,1})}\leq
\| rv_{0,\theta}\| _{L^{\infty}(\reals^3)}<\infty,
\end{equation}
where we have used
the hypothesis that $rv_{0,\theta}\in L^{\infty}(\reals^3)$, the
point-wise bound in Proposition \ref{chaeandleeprop}, and the fact that
$1<\sqrt{x_1^2+x_2^2}<4$. Let
\begin{equation}\Omegabar_+(x,t)=   \left  \{\begin{array}{cc}
 \Omega(x,t)+\Lambda & \Omega(x,t)\geq0, \\
 \Lambda & \Omega(x,t)<0. \\
\end{array}   \right  .\end{equation}
Note that $\Omegabar_+\geq\Lambda$ and all derivatives of
$\Omegabar_+$ on the set where $\Omega(x,t)<0$ are equal to zero.
This function is also Lipschitz and $\Omega$ is smooth by assumption.
At interfaces boundary terms upon integration by parts will cancel
and so we can make sense of the calculations below. Direct
computation yields, for $q>1$, that
\begin{equation}
\label{Omegaq}\Delta
\Omegabar_+^q -(b\cdot \nabla) \Omegabar_+^q
+\frac{2}{r}\partial_r \Omegabar_+^q-\partial _t
\Omegabar_+^q=-\frac{q\Omegabar_+^{q-1}}{r^2}\frac{\partial
v_{\theta}^2}{\partial z}+q(q-1)\Omegabar_+^{q-2}|\nabla
\Omegabar_+|^2.
\end{equation}
 Let $\frac{5}{8}\leq
\sigma_2<\sigma_1 \leq 1$. Define
\be
\al
\lab{defPsigma}
 P(\sigma_i)
&
 = \{(r,\theta,z) :
(5-4\sigma_i)<r<4\sigma_i,\ 0\leq\theta\leq 2\pi,\
|z|<4\sigma_i\}\times (-\sigma_i^2,0)\\
&=C(\sigma_i)
\times(-\sigma_i^2,0),
\eal
\ee
 for $i=1,2$.   Here for convenience
denote the space portion, which is a hollowed out cylinder, as
$C(\sigma_i)$ .
Choose $\psi=\phi(y)\eta(s)$ to be a refined cut-off function
satisfying
\begin{align*}
\supp\phi \subset C(\sigma_1)&;\ \phi(y)=1  \textrm{ for all }
y\in C(\sigma_2);\ 0\leq \phi\leq 1;\\& \frac{|\nabla \phi|}{\phi
^{\delta}}\leq\frac{c_1}{\sigma_1-\sigma_2}\textrm{ for } \delta
\in (0,1) \,  \text{to be chosen later in the proof};
\\
\supp \eta \subset (-\sigma_1^2,0];&\
\eta(s)=1,\textrm{ for all } s\in [-\sigma_2^2,0];\ 0\leq \eta \leq 1\\
&\ |\eta '|\leq
\frac{c_2}{(\sigma_1-\sigma_2)^2}.
\end{align*}
 Let
$f=\Omegabar_+^q$ and use $f\psi ^2$ as a test function in
(\ref{Omegaq}) to get
\begin{align*} \int_{P(\sigma_1)}(&\Delta f - (b\cdot\nabla)f -\partial _s f
+\frac{2}{r}\partial
_rf)f \psi ^2 dyds\\
&=
\int_{P(\sigma_1)}q(q-1)\Omegabar_+^{q-2} |\nabla \Omegabar_+|^2f \psi
^2 dyds-\int_{P (\sigma_1)}\frac{q\Omegabar_+^{q-1}}{r^2}\frac{\partial
v_{\theta}^2}{\partial z}f\psi^2dyds\\
\ &=q(q-1)\int_{P(\sigma_1)}\Omegabar_+ ^{-2}|\nabla
\Omegabar_+|^2f^2\psi ^2dyds-\int_{P
(\sigma_1)}\frac{q\Omegabar_+^{2q-1}}{r^2}\frac{\partial
v_{\theta}^2}{\partial z}\psi^2dyds\\
\ &\geq
-\int_{P(\sigma_1)}\frac{q\Omegabar_+^{2q-1}}{r^2}\frac{\partial
v_{\theta}^2}{\partial z}\psi^2dyds.
\end{align*}

Integration by parts on the first term implies that
\begin{align*}
&\int_{P(\sigma_1)}\nabla (f \psi ^2)\nabla f dyds\\
&\leq \int_{P(\sigma_1)}  \Big (-b\cdot\nabla f(f \psi ^2)-\partial
_s f(f \psi ^2) +\frac{2}{r}\partial _rf(f \psi ^2)
+\frac{q\Omegabar_+^{2q-1}}{r^2}\frac{\partial
v_{\theta}^2}{\partial z}\psi^2   \Big  )dyds.
\end{align*}
 A manipulation
using the product rule shows that
\begin{align*}
\int_{P(\sigma_1)}\nabla(f \psi ^2)\nabla f dyds
=\int_{P(\sigma_1)}
\left  (|\nabla(f \psi)|^2-|\nabla
\psi|^2f^2   \right  )dyds.
\end{align*}
Thus,
\begin{align*}
& \int_{P(\sigma_1)}|\nabla(f  \psi)|^2dyds
\leq \int_{P(\sigma_1)}  \Big (-b\cdot\nabla f(f \psi ^2)-\partial
_s f(f \psi ^2) +\frac{2}{r}\partial _rf(f \psi ^2)\\
&\hspace{5cm}
+\frac{q\Omegabar_+^{2q-1}}{r^2}\frac{\partial
v_{\theta}^2}{\partial
z}\psi^2+|\nabla\psi|^2f^2  \Big )dyds.
\end{align*}
 Integration by parts
on the term involving the time derivative yields
\begin{align*}
& \int_{P(\sigma_1)}-(\partial _s f)f \psi ^2dyds
 =-\frac{1}{2}\int_{P(\sigma_1)}\partial_s(f^2)\psi^2  dyds\\
&=-\tfrac{1}{2} \Big (\int_{C(\sigma_1)}
 \!  \! f^2\psi^2(y,0)dy
-\int_{C(\sigma_1)}
 \!  \!
f^2\psi^2(y,-\sigma_1
^2)dy  \Big  )
+\tfrac{1}{2}
\int_{P(\sigma_1)}
 \!  \!
\partial_s(\psi^2)f^2dyds.
\end{align*}
 Our
cut-off functions provide  $\psi^2=   (\phi\eta   )^2,\
\eta(0)=1,\
\eta(-\sigma_1 ^2)=0$, and $0\leq\phi\leq1$.  Thus,
\begin{align*}\int_{P(\sigma_1)}-(\partial _s f)f \psi ^2dyds&
=-\frac{1}{2}\int_{C(\sigma_1)}f^2(y,0)\phi^2(y)dy+
\int_{P(\sigma_1)}\phi^2(\eta\partial_s\eta)f^2dyds\\
&\leq-\frac{1}{2}\int_{C(\sigma_1)}f^2(y,0)\phi^2(y)dy
+\int_{P(\sigma_1)}(\eta\partial_s\eta) f^2dyds,
\end{align*}
and so,
\begin{equation}
\label{T1T2T3T4}
\begin{aligned}
\int_{P(\sigma_1)}&|\nabla(f
\psi)|^2dyds+\frac{1}{2}\int_{C(\sigma_1)}
f^2(y,0)\phi^2(y)dy\\&\leq\int_{P(\sigma_1)}-b\cdot
\nabla
f(f\psi^2)dyds
+\int_{P(\sigma_1)}(\eta\partial_s\eta+
|\nabla\psi|^2)f^2dyds\\&\hspace{.5cm}+\int_{P
(\sigma_1)}\frac{2}{r}\partial _r
f(f\psi^2)dyds+\int_{P(\sigma_1)}
\frac{q\Omegabar_+^{2q-1}}{r^2}\frac{\partial v_{\theta}^2}{\partial
z}\psi^2dyds
\\&:=T_1+T_2+T_3+T_4.
\end{aligned}
\end{equation}

We will bound each of the terms on the right hand side in the next few steps.

\textbf{Step 2:}
In this step we find an upper bound for  $T_1$, following an idea in \cite{Z}, where a parabolic
equation with a similar drift term is explored.
The new input is that we can exploit the fact that in the
space time domains of concern, all three dimensional integrals are equivalent to
two dimensional ones.
Therefore we can apply the 2 dimensional Sobolev inequality, which allows us to make gains.
Another idea is to derive an a priori bound for the mean oscillation of  the angular stream function
by exploiting the structure of the equations.

Since $\Div b=0$,
\begin{align*} T_1&=\int_{P(\sigma_1)}-b\cdot(\nabla f)(f
\psi^2)dyds\\
&=\frac{1}{2}\int_{P(\sigma_1)}-b\psi^2\cdot\nabla
(f^2)dyds=\frac{1}{2}\int_{P(\sigma_1)}\Div(b
\psi^2)f^2dyds\\
&=\frac{1}{2}\int_{P(\sigma_1)}\Div b (\psi
f)^2dyds+\frac{1}{2}\int_{P(\sigma_1)} b\cdot\nabla (
\psi  ^2)f^2dyds\\
&= \int_{P(\sigma_1)}b\cdot(\nabla\psi)\psi f^2dyds.
\end{align*}

In the following, we will often carry out computations on
 the following domains with 2 spatial dimension
\be
\lab{PbarCbar}
\al
\b C(\sigma_1) &= \{ (r, z) \, |  (r, \theta, z) \in C(\sigma_1) \}, \\
\b P(\sigma_1) &=  \{ (r, z, s) \, |  (r, \theta, z, s) \in P(\sigma_1) \}.
\eal
\ee We will also use the following notations.
\[
\b y = (r, z), \quad  d\b y= dr dz, \quad \text{if} \quad dy = r dr dz d\theta.
\]

Let $L_\theta$ be the angular component of the stream function in
the cylindrical coordinates. It is well known that
\be
\lab{vandL}
v_r = - \partial_z L_\theta, \quad v_z = \frac{1}{r} \partial_r (
r L_\theta).
\ee  Note
$L_\theta$ is also axially symmetric.

Let  $a=a(t)$
be a function of time only, which will be chosen later.
Using integration by parts and the divergence free property of
$b=v_r e_r + v_z e_z$, we have
\begin{eqnarray} \nonumber
  T_1 &=&\int_{P(\sigma_1)}b \cdot(\nabla \psi)(\psi f^2)dyds \\                              \nn
   &=&\int_{P(\sigma_1)} ( v_z \partial_z \psi + v_r \partial_r \psi) (\psi f^2)dyds \\               \nn
   &=&2 \pi \int_{\bar{P}(\sigma_1)} \partial_r (r L_\theta - a) \, (\partial_z \psi) (\psi f^2)dr dz ds
   -2 \pi \int_{\bar{P}(\sigma_1)} \partial_z (r L_\theta - a) (\partial_r \psi) (\psi f^2))dr dz ds\\                \nn
   &=& -2 \pi \int_{\bar{P}(\sigma_1)}  (r L_\theta - a) \, \partial_r [(\partial_z \psi) (\psi f^2)]dr dz ds
   +2 \pi \int_{\bar{P}(\sigma_1)}
 (r L_\theta - a) \partial_z[ (\partial_r \psi) (\psi f^2))] dr dz ds. \\                \nn
\end{eqnarray}
From here, a routine calculation shows
\begin{eqnarray} \label{8.1}
T_1&\le &\frac{1}{8}\int_{P(\sigma_1)}|\nabla (\psi f)|^2dyds  \\
\nn &&+C \frac{\sup\limits_{t\in (-\sigma^2_1,0)} \left(\|r
L_\theta- a(t)\|^2_{L^\infty(\bar{C}(\sigma_1))} +
1\right)}{(\sigma_1-\sigma_2)^2}
\int_{\bar{P}(\sigma_1)}f^2d\bar{y}ds.
\end{eqnarray} Here we have used the fact that $r$ is comparable to $1$ in the region
$P(\sigma_1)$ and therefore, the 2 dimensional volume element $dr dz$ is comparable to
the three dimensional one $r dr dz d\theta$ for axially symmetric functions.

Choose a  2 dimensional cut-off function $\phi= \phi(r, z)\in C^\infty_0(\reals^2)$
such that $\phi=1$ in $\bar{C}(\sigma_1)$, $supp \, \phi \in
\bar{C}(9\sigma_1/8)$ and  $0 \le \phi \le 1$, and $|\bar{\nabla}
\phi| + |\Delta_2 \phi| \le C$.  Observe
that the region $\bar{C}(9\sigma_1/8)$ is  at least $1/2$ unit away from the
$z$ axis. Here and later in this section, $\bar{\nabla} =
(\partial_r,
\partial_z)$ is the 2 dimensional gradient, and $\Delta_2 =
\partial^2_r + \partial^2_z$ is the $2$ dimensional Laplacian with
respect to the $r$ and $z$ variables. In the $2$ dimensional space
of $(r, z)$, according to Hou-Li  \cite{HL2} Appendix,
\begin{eqnarray} \nonumber
\|r L_\theta-a(t)\|_{L^\infty(\bar{C}(\sigma_1))}&\leq &\|(r
L_\theta-a(t))\phi\|_{L^\infty(\bar{C}(9\sigma_1/8))}  \\    \nn
     &\leq& C  \left (\|\bar{\nabla}((r L_\theta-a(t))\phi)\|_{L^2(\reals^2)}+
\|(r L_\theta-a(t))\phi\|_{L^2(\reals^2)}+1 \right) \times  \\
\nn
    && \left[\log(\|\Delta_2((r L_\theta-a(t))\phi)
     \|_{L^2(\reals^2)}+\|(r L_\theta-a(t))\phi\|_{L^2(\reals^2)}+e)\right]^\frac{1}{2} \\  \nn
&\leq&C \left( \|\phi \bar{\nabla}(r L_\theta-a(t))\|_{L^2(\reals^2)}+
\|(r L_\theta-a(t)) \|_{L^2(\bar{C}(9\sigma_1/8))}+1 \right) \times
\\  \nn
    && \left[\log(\|\Delta_2((r L_\theta-a(t))\phi)
     \|_{L^2(\reals^2)}+\| (r L_\theta-a(t)) \|_{L^2(\bar{C}(9\sigma_1/8))}+e)\right]^\frac{1}{2}.
\end{eqnarray}

Choose $a(t)$ be the average of $r L_\theta(\cdot, t)$ on
$\bar{C}(9\sigma_1/8)$ under the $2$ dimensional volume element $d r
dz$.  Using the $2$ dimension Poincar\'{e} inequality, we deduce
\begin{eqnarray} \label{8.2}
&&\|r L_\theta-a(t)\|_{L^\infty(\bar{C}(\sigma_1))} \leq  C
 (\|\bar{\nabla}  (r L_\theta) \|_{L^2(\bar{C}(9\sigma_1/8))}+1)\times  \\
\nn
    && \qquad \left[\log(\|\Delta_2((r L_\theta -a(t))\phi)
     \|_{L^2(\reals^2)}+ C \|\bar{\nabla} (r L_\theta) \|_{L^2(\bar{C}(9\sigma_1/8))}+e)\right]^\frac{1}{2}
\end{eqnarray}

Note that at any point $(x, t)$ we have, \be \lab{drl=rb}
|\bar{\nabla} (r L_\theta) |^2 = \left| \partial_r ( r L_\theta)
\right|^2 + r^2 \left| \partial_z L_\theta \right|^2 = r^2 \left[
\left|\frac{1}{r} \partial_r ( r L_\theta) \right|^2 + \left|
\partial_z L_\theta \right|^2 \right] = r^2 |b |^2.
\ee

Also
\[
\al \Delta_2 &((r L_\theta-a(t))\phi)  = \Delta_2(r L_\theta )
\phi + 2 \bar{\nabla} (r L_\theta) \cdot \bar{\nabla} \phi + (r
L_\theta-a(t)) \Delta_2 \phi\\
&= \phi (\partial^2_r+\partial^2_z) (r L_\theta) + 2 \bar{\nabla}
(r L_\theta) \cdot \bar{\nabla} \phi + (r L_\theta-a(t)) \Delta_2
\phi\\
&=\phi r \left( \partial^2_r L_\theta +\partial^2_z L_\theta +
\frac{2}{r}
\partial_r L_\theta \right) + 2 \bar{\nabla} (r L_\theta) \cdot \bar{\nabla} \phi
+ (r L_\theta-a(t)) \Delta_2 \phi\\
&=\phi r \left( \partial^2_r L_\theta +\partial^2_z L_\theta +
\frac{1}{r}
\partial_r L_\theta - \frac{1}{r^2} L_\theta \right) +
 \phi \left( \partial_r L_\theta + \frac{1}{r} L_\theta \right)+ 2 \bar{\nabla} (r L_\theta) \cdot \bar{\nabla} \phi
+ (r L_\theta-a(t)) \Delta_2 \phi\\
&=-\phi r w_\theta + \phi v_z + 2 \bar{\nabla} (r L_\theta) \cdot
\bar{\nabla} \phi + (r L_\theta-a(t)) \Delta_2 \phi. \eal
\] Therefore, together with (\ref{drl=rb}), we have,
\[
|\Delta_2 ((r L_\theta-a(t))\phi)| \le |r w_\theta + v_z| + 2 r
|b| + C |(r L_\theta-a(t))|.
\] Using the 2 dimensional Poincar\'e inequality again, we find
that
\[
\al
 \|\Delta_2&((r L_\theta -a(t))\phi)
     \|_{L^2(\reals^2)} \\
 &\le C \|w_\theta \|_{L^2(\bar{C}(9\sigma_1/8))} +
     C \| b \|_{L^2(\bar{C}(9\sigma_1/8))} +  C \|(r
     L_\theta-a(t))\|_{L^2(\bar{C}(9\sigma_1/8))}\\
 &\le C \|w_\theta \|_{L^2(\bar{C}(9\sigma_1/8))} +
     C \| b \|_{L^2(\bar{C}(9\sigma_1/8))} +  C \| \bar{\nabla}  (r
     L_\theta)\|_{L^2(\bar{C}(9\sigma_1/8))}\\
&\le C \|w_\theta \|_{L^2(\bar{C}(9\sigma_1/8))} +
     2 C \| b \|_{L^2(\bar{C}(9\sigma_1/8))}. \eal
\]

 Then from \eqref{8.2}, we get
\be
\lab{oscL}
\al
&\|r L_\theta- a(t) \|_{L^\infty(\bar{C}(\sigma_1))}\\
&\leq C
(\|v(\cdot,t)\|_{L^2(\bar{C}(9\sigma_1/8))}+1)
    \left[\log(C \|w_\theta(\cdot,t)\|_{L^2(\bar{C}(9\sigma_1/8))}
     +C \|v(\cdot,t)\|_{L^2(\bar{C}(9\sigma_1/8))}+e)\right]^\frac{1}{2}.
\eal
\ee
Then by \eqref{8.1}, we have
\[
\al
&T_1\le \frac{1}{8}\int_{P(\sigma_1)}|\nabla (\psi f)|^2dyds  \\
&+ C \frac{\sup\limits_{t\in
(-\sigma^2_1,0)}\left[ (\|v(\cdot,t)\|_{L^2(\bar{C}(9\sigma_1/8))}+1)^2
    \log(\|w_\theta(\cdot,t)\|_{L^2(\bar{C}(9\sigma_1/8))}
     +\|v(\cdot,t)\|_{L^2(\bar{C}(9\sigma_1/8))}+e) \right]}{(\sigma_1-\sigma_2)^2}
\int_{\bar{P}(\sigma_1)}f^2d\bar{y}ds.
\eal
\]
Using the notation
\be
\lab{Kbar}
\bar K=\bar{K}(v, w) \equiv \sup\limits_{t\in
(-\sigma^2_1,0)}\left[ (\|v(\cdot,t)\|_{L^2(\bar{C}(9\sigma_1/8))}+1)
    \log^{1/2}(\|w_\theta(\cdot,t)\|_{L^2(\bar{C}(9\sigma_1/8))}+\|v(\cdot,t)\|_{L^2(\bar{C}(9\sigma_1/8))}+e) \right]
\ee we can write the last inequality as
\be
\lab{t1<<}
T_1 \le \frac{1}{8}\int_{P(\sigma_1)}|\nabla (\psi
f)|^2dyds+C \frac{\bar K^2(v, w)}{(\sigma_1-\sigma_2)^2}
\int_{\bar{P}(\sigma_1)}f^2d\bar{y}ds.
\ee

 \textbf{Step 3:} The term
$T_2$  is treated routinely.  We use
  \begin{align*}
T_2=\int_{P(\sigma_1)}(\eta\partial_s\eta+|\nabla\psi|^2)f^2dyds,
  \end{align*}
and properties of the cutoff,
  \begin{align*}
|\nabla \psi |^2=|\eta \nabla \phi |^2\leq  \Big ( \frac{|\nabla
\phi |}{\phi ^{\delta }}  \Big ) ^2\leq
\frac{c_1^2}{(\sigma_1-\sigma_2)^2},
  \end{align*}
and
  \begin{align*}
|\eta
\partial _s \eta |\leq |\partial _s \eta |
\leq \frac{c_2}{(\sigma_1-\sigma_2)^2},
\end{align*}
to get
\begin
{equation}\label{T2}
|T_2|\leq
\frac{C}{(\sigma_1-\sigma_2)^2}\int_{P(\sigma_1)}
f^2 dyds \le \frac{C}{(\sigma_1-\sigma_2)^2}\int_{\b P(\sigma_1)}
f^2 d\b yds.
\end{equation}

\textbf{Step 4:} As we deal with
$T_3=\int_{P(\sigma_1)}\frac{2}{r}\partial _rf(f\psi^2)dyds$, we
note we are assuming the integration takes place away from the
singularity set of the solution to the axially symmetric Navier
Stokes equations and away from the z-axis in general. Thus, all
functions are bounded and smooth and $r$ varies between two
positive constants.  We also utilize the cylindrical coordinates
of the axially symmetric case, and integration by parts:
  \begin{align*}
 T_3
&=\int_{P(\sigma_1)}\frac{2}{r}\partial
_rf(f\psi^2)dyds
=\int_{P(\sigma_1)}\frac{1}{r}\partial_r(f^2)\psi^2rdrd\theta
dzds
\\&
=\int_{P(\sigma_1)}\partial_r(f^2)\psi^2drd\theta
dzds
=-\int_{P(\sigma_1)}\partial_r(\psi^2)f^2drd\theta
dzds
\\&
=-\int_{P(\sigma_1)}\frac{2}{r}\partial_r\psi (\psi
f^2)rdrd\theta dzds
 =-\int_{P(\sigma_1)}\frac{2}{r}\partial_r(\psi)(\psi f^2)dyds\\
&=-\int_{P(\sigma_1)}\frac{2}{r}\overrightarrow{e_r}\cdot\nabla\psi(\psi
f^2)dyds.
\end{align*}
 The Cauchy-Schwartz inequality then implies
\begin{align*}
|T_3|\leq\int_{P(\sigma_1)}\frac{2}{r}|\nabla\psi|\psi f^2dyds.
\end{align*} This yields

\begin{equation}\label{T3}
|T_3| \leq
\frac{C}{(\sigma_1-\sigma_2)} \int_{P(\sigma_1)}f^2dyds \le \frac{C}{(\sigma_1-\sigma_2)} \int_{\b P(\sigma_1)}f^2 d\b yds .
\end{equation}

 \textbf{Step 5:}
Lastly, we work on the inhomogeneous
term of (\ref{Omega/omegatheta}), that is,
$\displaystyle\frac{2v_{\theta}}{r^2}\frac{\partial
v_{\theta}}{\partial z}$, which produced the term $T_4$. Recall
\begin{align*}\Lambda=\| v_{\theta}\| _{L^{\infty}(P_{1,4,1})}\leq
\| rv_{0,\theta}\| _{L^{\infty}(\reals^3)}<\infty,\end{align*}
 and that
$\Omegabar_+=
\begin{cases} \Omega+\Lambda & \Omega \geq
0\\\Lambda & \Omega<0\end{cases}      $, thus
$\Omegabar_+\geq\Lambda$. Also, we have let $f=\Omegabar_+ ^q$.
Using integration by parts yields
\begin{align*}
T_4=&\int_{P(\sigma_1)}\frac{q\Omegabar_+^{2q-1}}{r^2}\frac{\partial
v_{\theta}^2}{\partial z}\psi^2dyds\\
=&-\int_{P(\sigma_1)}\frac{\partial}{\partial z}
    \Big   (\frac{\Omegabar_+^{2q} \psi^2}{\Omegabar_+}      \Big   )
\frac{q}{r^2}
v_{\theta}^2  dyds\\
=& -\int_{P(\sigma_1)}\frac{\partial}{\partial z} (f\psi)^2
\frac{1}{\Omegabar_+}  \frac{q}{r^2} v_{\theta}^2  dyds +
\int_{P(\sigma_1)} (\Omegabar_+^q \psi)^2 \frac{1}{\Omegabar_+^2}
\frac{\partial \Omegabar_+}{\partial z} \frac{q}{r^2} v_{\theta}^2
dyds \\
=& -\int_{P(\sigma_1)}\frac{\partial}{\partial z} (f\psi)^2
\frac{1}{\Omegabar_+}  \frac{q}{r^2} v_{\theta}^2  dyds \\
&\hspace{3cm}+
\frac{1}{2} \int_{P(\sigma_1)} \frac{1}{\Omegabar_+}
    \Big   [\frac{\partial (\Omegabar_+^{2q}\psi^2)}{\partial z} -
\Omegabar_+^{2q} \frac{\partial \psi^2}{\partial z}     \Big   ]
\frac{1}{r^2} v_{\theta}^2 dyds
\\
=& -\int_{P(\sigma_1)}\frac{\partial}{\partial z} (f\psi)^2
\frac{1}{\Omegabar_+}  \frac{q-(1/2)}{r^2} v_{\theta}^2  dyds -
\frac{1}{2} \int_{P(\sigma_1)} \frac{1}{\Omegabar_+}
\Omegabar_+^{2q} \frac{\partial \psi^2}{\partial z} \frac{1}{r^2}
v_{\theta}^2 dyds.
\end{align*}
Considering that  $\frac{|v_{\theta}|}{\Lambda}\leq 1$, utilizing
$\Lambda\leq\Omegabar_+$, and $r=\sqrt{y_1^2+y_2^2}\geq 1$ for all
$y\in P(\sigma_1)$, we continue by fixing $\epsilon_3>0$. Apply
Young's inequality with exponents both being $2$ to get
\[
\al
&|T_4|
 \leq\int_{P(\sigma_1)}2q|v_{\theta}\| f|\psi
  \Big   |\frac{\partial(f\psi)}{\partial
z}     \Big   |dyds+ \frac{c_3}{\sigma_1-\sigma_2} \int_{P(\sigma_1)}
f^{2}
|v_{\theta}| dyds\\
\notag
&\leq \int_{P(\sigma_1)}
    \Big   |\frac{2q\Lambda}{(2\epsilon_3)^{\frac{1}{2}}}f\psi     \Big   |
\times
    \Big   |(2\epsilon_3)^{\frac{1}{2}}\frac{\partial(f\psi)}{\partial
z}     \Big   | dyds + \frac{c_3 \Lambda}{\sigma_1-\sigma_2}
\int_{P(\sigma_1)} f^{2}
 dyds\\
\notag
&\leq \frac{c_{12}\Lambda^2q^2}{\epsilon_3}\int_{P(\sigma_1)}f^2dyds
+\epsilon_3\int_{P(\sigma_1)}|\nabla(f\psi)|^2dyds+ \frac{c_3
\Lambda}{\sigma_1-\sigma_2} \int_{P(\sigma_1)} f^{2}
 dyds.
\eal
\]Thus
\be
\lab{t4<}
|T_4| \le \frac{1}{4} \int_{ P(\sigma_1)}|\nabla(f\psi)|^2 dyds + C\left[\Lambda^2 q^2
+  \frac{\Lambda}{\sigma_1-\sigma_2} \right]  \int_{\b P(\sigma_1)}f^2 d\b yds.
\ee
\medskip

\textbf{Step 6: $L^2-L^{\infty}$ Estimate:} An $L^2-L^{\infty}$
bound is derived using Moser's iteration. Recall inequality
(\ref{T1T2T3T4}) from Step 1  and substitute the estimates for $T_1, T_2, T_3,
T_4$ in
(\ref{t1<<}), (\ref{T2}),   (\ref{T3}),  (\ref{t4<}) respectively,
we obtain
\begin{align*}
&\int_{P(\sigma_1)}|\nabla(f\psi)|^2dyds+
\frac{1}{2}\int_{C(\sigma_1)}f^2(y,0)\phi^2(y)dy\\
&\leq  \frac{3}{4} \int_{P(\sigma_1)}|\nabla(f\psi)|^2 dyds
+ C\left[\frac{C \b K^2(v, w) }{(\sigma_1-\sigma_2)^2}  + 1+ \Lambda^2 q^2
+  \frac{1}{(\sigma_1-\sigma_2)^2} \right]  \int_{\b P(\sigma_1)}f^2 d\b yds.
\end{align*}  Here $\b K=\bar K(v, w)$ is defined in (\ref{Kbar}).

Consequently,
\begin{align}
\label{energyestimate}
\int_{\b P(\sigma_1)}|\b \nabla(f\psi)|^2&d\b yds
+\int_{\b C(\sigma_1)}f^2(y,0)\phi^2(y)d\b y\\
\notag
&\leq\frac{C q^2}{(\sigma_1-\sigma_2)^2}
   \Big   (\b K^2(v, w) (C_{1,4})+\Lambda^2+1     \Big   )\int_{\b P
(\sigma_1)}f^2d\b yds.
\end{align}
The last inequality follows since $q>1$ and
$0<\sigma_1-\sigma_2<1$.

\smallskip

Next we will iterate the above energy type estimate.
We will apply  refined interpolation and embedding inequalities involving
BMO functions.
Applying Lemma 1 in Section 2 of \cite{KT}, we know that
\[
\Vert (f \phi)^2 \Vert_{L^2(\b C(\sigma_1))} \le C \Vert f \phi \Vert_{L^2(\b C(\sigma_1))} \Vert f \phi\Vert_{BMO(\b C(\sigma_1))}.
\]By the 2 dimensional Poincar\'e inequality, we also have
\[
\Vert f \phi \Vert_{BMO(\b C(\sigma_1))} \le C \Vert \b \nabla (f \phi) \Vert_{L^2}.
\] These two inequalities imply that
\[
\Vert f \phi \Vert_{L^4(\b C(\sigma_1))} \le C \Vert f \phi \Vert^{1/2}_{L^2(\b C(\sigma_1))} \Vert\b \nabla( f \phi) \Vert^{1/2}_{L^2(\b C(\sigma_1))}.
\]Consequently
\[
\al
\int_{-\sigma_1^2}^0\int_{\reals^2}(&f\psi)^4 d\b yds\\&\leq
C \sup_{-\sigma_1^2
\leq s\leq 0}     \Big   (\int_{\reals^2}(f\psi)^2 d\b y      \Big
) \int_{-\sigma_1^2}^0\int_{\reals^2}|
\b \nabla(f\psi)|^2d\b yds.
\eal
\]Substituting (\ref{energyestimate}) into the right hand side of the above inequality, we
deduce
\[
\al
\int_{\b P(\sigma_2)} f^4 d \b y ds \le C \left[ \frac{q^2}{(\sigma_1-\sigma_2)^2}
 \left[\bar K^2
  + \Lambda^2
+  1 \right]  \int_{\b P(\sigma_1)}f^2 d\b yds \right]^2.
\eal
\] Here we have written $\bar K = \bar K(v, w)$ for simplicity.
This shows, since $f = \b \Omega^q_+$ by definition, that
\be
\lab{2MoserIt}
\al
\int_{\b P(\sigma_2)} \b \Omega^{4q}_+ d \b y ds
\le C \left[ \frac{q^2}{(\sigma_1-\sigma_2)^2}
 \left[\bar K^2
  + \Lambda^2
+  1 \right]  \int_{\b P(\sigma_1)} \b \Omega^{2q}_+ d\b yds \right]^2.
\eal
\ee

For $i=0, 1, 2, ...$,  in (\ref{2MoserIt}), we take $q=2^i$, replace
$\sigma_1$ by $\sigma_i=1-\Sigma^i_{j=1} 2^{-j-2}$ and $\sigma_2$ by
$\sigma_{i+1} =1-\Sigma^{i+1}_{j=1} 2^{-j-2}$.
Here we set $\Sigma^i_{j=1} 2^{-j-2}=0$ when $i=0$. Then
\begin{equation}
\label{mosergammai}
\left( \int_{\b P(\sigma_{i+1})}\Omegabar_+^{2^{i+2}}d\b yds \right)^{1/2}
 \leq
c_1     c^{i+1}_2 2^{2 i}      \Big
(\b K^2+\Lambda^2+1     \Big   ) \int_{\b P(\sigma_i)}
\Omegabar_+^{2^{i+1}}d\b yds.
\end{equation}

After iteration, we
 obtain
\begin{align*}
   \Big   (&\int_{\b P(\sigma_{i+1})}
 \Omegabar_+^{2^{i+2}}d\b yds     \Big   )^{\frac{1}{2^{i+1}}}\\
&\leq c_1^{\sum\frac{1}{2^{j}}}c_2^{\sum \frac{j+1}{
2^{j-1}}} 2^{2\sum \frac{j-1}{2^{j-1}}}
  \big
(\b K^2+
\Lambda^2+1     \big   )^{\sum\frac{1}{2^{j-1}}}
\int_{\b P_{1,4,1}}\Omegabar_+^2 d\b yds.
\end{align*}
 Note the sums in the exponents are all from $j=1$ to
$j=i+1$. Let $i     \rightarrow \infty$. All the exponential series
converge and we deduce.
\[
\sup_{\b P_{2, 3, 3/4}} \b \Omega^2_+
 \le C \left(\bar K^2
  + \Lambda^2
+  1 \right)^2   \int_{\b P_{1, 4, 1}} \b \Omega^2_+ d\b yds.
\] Repeating the argument on
$\Omegabar_-=
\left \{\begin{array}{cc}-\Omega+\Lambda &
\Omega\leq0\\
\Lambda&\Omega>0\end{array} \right.$,  we see that
\[
\sup_{\b P_{2, 3, 3/4}} \b \Omega^2_- \le
C \left(\b K^2 +\Lambda^2
+  1 \right)^2   \int_{\b P_{1, 4, 1}} \b \Omega^2_- d\b yds.
\] Thus
\[
\sup_{\b P_{2, 3, 3/4}}  \Omega^2
\le C  \left(\b K^2 +\Lambda^2
+  1 \right)^2   \int_{\b P_{1, 4, 1}} ( \Omega^2 + \Lambda^2)  d\b yds.
\]
In the region $P_{1, 4, 1}$,
 the quantities $\Omega=\omega_\theta/r$ and $\omega_\theta$ are equivalent.
 Hence, the above implies
\[
\sup_{\b P_{2, 3, 3/4}}  \omega_\theta^2 \le
C \left(\b K^2 +\Lambda^2
+  1 \right)^2  \int_{\b P_{1, 4, 1}} ( \omega_\theta^2 + \Lambda^2)  d\b yds.
\]

\textbf{Re-scaling:} We now recall
that we omitted the ``tildes" in the notation in the above
computations.  So what has actually been proven thus far is
\be
\lab{mvitO}
\sup_{(\~{x},\~{t})\in P_{2,3,\frac{3}{4}}}\~{\omega}^2_{\theta}
(\~{x},\~{t})\leq C    \Big   (\bar K +\~{\Lambda}+1     \Big   )^4      \Big
(\int_{P_{1,4,1}}\~{\omega}^2_{\theta}(\~{x},\~{t})
d\~{x}d\~{t}+\~{\Lambda}^2     \Big   ),
\ee
where $\~{x}=\frac{x}{k},\
\~{t}=\frac{t}{k^2}$ and from (\ref{Kbar}) with $\sigma_1=1$ that,
\[
\al
&\bar{K} = \bar{K}(\~ v, \~ w) \\
&\equiv \sup\limits_{\~ t \in
(-1, 0)}\left[ (\|\~v(\cdot, \~ t)\|_{L^2(\bar{C}(9/8))}+1)
    \log^{1/2}(\|\~ w_\theta(\cdot, \~ t)\|_{L^2(\bar{C}(9/8))}+
\|\~ v(\cdot, \~t)\|_{L^2(\bar{C}(9/8))}+e) \right]\\
& \le \sup\limits_{\~ t \in
(-1, 0)}\left[ (\|\~v(\cdot, \~ t)\|_{L^2(\bar{C}(0.5, 5))}+1)
    \log^{1/2}(\|\~ w_\theta(\cdot, \~ t)\|_{L^2(\bar{C}(0.5, 5))}+
\|\~ v(\cdot, \~t)\|_{L^2(\bar{C}(0.5, 5))}+e) \right].
\eal
\]Here we have used the fact that $\bar{C}(9/8) \subset \bar{C}(0.5, 5)$ from
the definition in
(\ref{defPsigma}),  (\ref{R3tong}) and the 2 dimensional versions (\ref{PbarCbar}).
 We recall from the beginning of the section
 that
  \begin{align*}
\|\~{b}(\~{x},\~{t})\|_{L^{\infty}(-1,0;
L^2(C_{0.5, 5}))}=k^{-\frac{1}{2}} \|
b(x,t)\|_{L^{\infty}{(-k^2,0;L^2(C_{0.5 k, 5 k}))}},
\end{align*}
\[
\|\~{w}(\~{x},\~{t})\|_{L^{\infty}(-1,0;
L^2(C_{0.5, 5}))}=k^{\frac{1}{2}} \|
w(x,t)\|_{L^{\infty}{(-k^2,0;L^2(C_{0.5 k, 5 k}))}},
\]
 and
  \begin{align*}
\|\~{\omega}(\~{x},\~{t})\|_{L^2(P_{1,4,1})}=
k^{-\frac{1}{2}}\|\omega(x,t)\|_{L^2(P_ {k,4k,k})}.
\end{align*}  Hence
\[
\bar K \le \sup_{t \in [-k^2, 0]}
\left( \frac{1}{k^{1/2}}\|b\|_{L^2(C_{1k, 4k})} + 1 \right) \log^{1/2}
(k^{1/2} \| w_\theta(\cdot, t)\|_{L^2(C_{0.5k, 5k})}+ k^{-1/2}
\| v(\cdot, t)\|_{L^2(C(0.5k, 5k))}+e)
\]From \cite{BZ:1}, we have
\[
|\omega_\theta(x, t)|\le \frac{C}{k^5}, \quad  x \in C(0.5k, 5k)
\] where $C$ only depends on the initial value.
So, together with the energy inequality we have, for $k \in (0, 1]$,
\be
\lab{barK<}
\bar K \le C (k^{-1/2} \|b\|_{L^\infty(-k^2, 0;  L^2(C_{1k, 4k}))} +1) \log^{1/2} (\frac{1}{k} +e)
\ee where $C$ depends only on the initial value.

Also,   we note the control on $\~\Lambda$ is a scaling invariant
quantity, since, by Proposition \ref{chaeandleeprop}, we have
\begin{align*}
 \~{\Lambda}&=     \Big
(\sup_{P_{1,4,1}}|\~{v_{\theta}}(\~{x},\~{t})|     \Big   )
\leq
\Big   (\| \~{r}\~{v_{\theta}}(\~{x},-T)\|_{L^{\infty}(\reals^3)}
\Big ) \\
&=\|rv_{\theta}(x,-T)\|_{L^{\infty}(\reals^3)}
 =\|rv_{0,\theta}\|_{L^{\infty}(\reals^3)}.
\end{align*} Substituting this and (\ref{barK<}) into (\ref{mvitO}),
we obtain, for $k \in (0, 1)$ that
 \begin{align*}
\sup_{(x,t)\in P_{2k,3k,\frac{3k}{4}}}k^4\omega ^2_{\theta}(x,t)&\leq
C     \Big   (  (k^{-1/2} \|b\|_{L^\infty(-k^2, 0;  L^2(C_{1k, 4k}))} +1)
\log^{1/2} (\frac{1}{k} +e) +\|rv_{0,\theta}\|_{L^{\infty}
(\reals^3)}  +1   \Big   )^4\\
&\hspace{2cm}   \times  \Big
(\int_{P_{k,4k,k}}k^4\omega_{\theta}^2(x,t)\frac{1}{k^5}dxdt+\|
rv_{0,\theta}\|_{L^{\infty}(\reals^3)}^2     \Big   )\\ &\leq
\frac{C}{k^{3}}     \Big   (\|b\|_{L^{\infty}(-k^2,0;
L^2({C_{1k, 4k}}))} \log^{1/2} (\frac{1}{k} +e) +k^{1/2} \|rv_{0,\theta}\|_{L^{\infty} (\reals^3)}
+ k^{1/2} \Big   )^4 \\
&\hspace{3.5cm}   \times \Big
(\|\omega_{\theta}\|_{L^2(P_{k,4k,k})}^2+k\|rv_{0,\theta}\|_{L^{\infty}
(\reals^3)}^2     \Big   ).
\end{align*}
Therefore,
\be \lab{wtheta<3.5} \al
\|\omega_{\theta}(x,t)\|_{L^{\infty}(P_{2k,3k,\frac{3k}{4}})}
&\leq\frac{C}{k^{7/2}}     \Big   (\|b\|_{L^{\infty}
(-k^2,0;L^2(C_{1k, 4k}))} \log^{1/2} (\frac{1}{k} +e) +k^{1/2}
\|rv_{0,\theta}\|_{L^{\infty}(\reals^3)} + k^{1/2}
   \Big   )^{2}  \\
&\hspace{3cm}    \times \left
 (\|\omega_{\theta}\|_{L^2(P_{k,4k,k})}+\sqrt{k}\|rv_{0,
\theta}\|_{L^{\infty}(\reals^3)}     \right).
\eal
\ee
This proves Theorem \ref{ourtheorem}. \qed
\medskip

\section{Velocity Bound, proof of Theorem \ref{ourtheoremV} }

Based on the bound on $\omega_\theta$, now we prove the a priori
velocity bound. In this section we ignore the time variable.  Let $b
= v_re_r + v_ze_z$. Then taking the cylindrical curl of $b$ we get
    \begin{align*}
        \text{curl }b= \omega_\theta e_\theta
    \end{align*}
Taking the curl on both sides of this equation and using the
divergence free condition on $b$ to get $\text{curl (curl } b) =
-\Delta b$, we have
    \begin{align*}
        -\Delta b = \text{curl (}\omega_\theta e_\theta)
    \end{align*}
Choose a smooth cutoff function $\phi$ with support contained in the
ball $B_{2r_0}=B(x,2r_0)$, $0 \leq \phi \leq 1$ in $B_{2r_0}$,
$\phi\equiv1$ in $B_{\frac{r_0}{4}}$ and the following properties:
    \begin{align*}
        |\nabla\phi|\leq\frac{C}{r_0}\quad,\quad |\Delta\phi|\leq\frac{C}{r_0^2}
    \end{align*}
Then $\text{supp}(\nabla\phi)\subset B_{2r_0}\setminus
B_{\frac{r_0}{4}}$, and we compute
    \begin{align*}
        \Delta(\phi b) &= \Delta\phi b + 2\nabla\phi\cdot\nabla b + \phi\Delta b \\
    \Rightarrow \phi b &= \int_{\ballt}{\Gamma(x,y)\Big(\Delta\phi\, b + 2\nabla\phi\cdot\nabla b - \phi\text{ curl (}\omega_\theta e_\theta)\Big)}\,dy \\
        \phi b &= \int_{\ballt}{\Gamma(x,y)\Delta\phi\, b}\,dy + 2\int_{\ballt}{\Gamma(x,y)\nabla\phi\cdot\nabla b}\,dy - \int_{\ballt}{\Gamma(x,y)\phi\text{ curl (}\omega_\theta e_\theta)}\,dy \\
    \end{align*}
where $\Gamma(x,y) = \displaystyle\frac{c_0}{|x-y|}$ is the Green's
function.  After doing integration by parts, it is easy to see that, for all $p \ge 1$,
\be
\lab{btow}
\sup_{B_{r_0}(x)} |b| \le C r^{-3/p}_0 \Vert b \Vert_{L^p(B_{2 r_0}(x))}
+ C r_0 \sup_{B_{2 r_0}(x)} |\omega_\theta|.
\ee

Now we pick a point $x=(x_1, x_2, x_3)$ and let $r=|x'|$ be the distance from $x$ to the $z$ axis again.
We can assume $r \le 1/2$ since we are only concerned with the bound near $z$ axis.
Take $r_0= r^{1.5} |\ln r|^{-1/2}$.  Then we can find
$r/(2 r_0) = r^{-1/2} |\ln r|^{1/2}/2$ (round up to nearest integer) many balls,
which are disjoint, and which
are generated by rotating $B(x, r_0)$ around the $z$ axis, such that,
their union is contained in the torus around the curve
\[
\{ y=(y_1, y_2, y_3) \, | \, \sqrt{y^2_1+y^2_2}=r,    y_3=x_3  \},
\] with cross sections being 2 dimensional balls of radius $r$.
Since the function $b$ is axially symmetric, the integral of $| b |^2$ on each of the ball
is the same number. Therefore
\[
\Vert b \Vert^2_{L^2(B_{2 r_0}(x))} \le C r^{1/2}  |\ln r|^{-1/2} \Vert b \Vert^2_{L^2(\reals^3)}.
\]

Using this, we can now take $p=2$  in (\ref{btow}) and use the a
priori bound on $\omega_\theta$ in Theorem \ref{ourtheorem} to deduce
\[
\al
|b(x, t)| &\le C r_0^{-3/2}   r^{1/2}  |\ln r|^{-1/2} + C r_0  r^{-7/2} |\ln r |\\
&=C r^{-9/4} |\ln r|^{3/4} r^{1/4}  |\ln r|^{-1/4} + C r^{3/2} |\ln r|^{-1/2}
r^{-7/2} |\ln r |\\
&= C r^{-2} | \ln r|^{1/2}.
\eal
\]
This completes the proof of Theorem \ref{ourtheoremV}. \qed

\section{Appendix: a regularity criteria for $v_z$ only}

In this section, we present 2 short results on equation (\ref{eqasns}), which may be of
independent interest.

The first one is a  critical regularity condition on the $z$ component of the velocity only.
We note that a similar result under a little more assumption is given in Theorem 2.2 of the
interesting paper
\cite{JX2}.  This result gives a mathematical explanation of the folklore belief that singularity
can happen only if the vertical convection is high enough, in a tornado e.g.

\begin{proposition}
Let $v$ be a Leray-Hopf solution to (\ref{eqasns}) in $\reals^3 \times (0, \infty)$ such that
$r v_{0, \theta} \in L^\infty(\reals^3)$.
Suppose, for a given constant $C>0$,
and all $x \in \reals^3$ and $t \ge 0$,
\[
|v_z(x, t)| \le \frac{C}{r}.
\]Then $v$ is regular for all time.
\end{proposition}

\proof  Recall from (\ref{vandL}) the well known relation $ v_z = \frac{1}{r} \partial_r (
r L_\theta)$.  Hence
\[
\al
\left| |x'| L_\theta(x, t) \right|
&= \left| \int^{|x'|}_0 \partial_r ( r L_\theta) dr  \right| \le
 \int^{|x'|}_0 \left| r v_z \right|  dr \le C |x'|.
\eal
\]Here we just used the  assumption on $v_z$.  Therefore $L_\theta$ is a bounded function.
Then from the main result in \cite{LZ11}, we know that $v$ is regular for all time. \qed

The next result is an a priori bound for the angular stream function, which scales in the same
way as the energy estimate modulo a log term.

\begin{proposition}
Under the same assumption as Theorem \ref{ourtheoremV}, there exists a constant $C$,
depending only on the initial data, such that,
\[
|L_\theta(x, t)| \le \frac{C |\ln |x'| |^{1/2}}{|x'|^{1/2}}, \qquad |x'| \le \min \{1/2, R \}.
\]
\end{proposition}
\proof
Recall from
(\ref{oscL}) the a priori bound for  the scaled $L_\theta$ in a dyadic cube at time $t$:
\[
\al
&\|r L_\theta- a(t) \|_{L^\infty(\bar{C}(\sigma_1))}\\
&\leq C
(\|v(\cdot,t)\|_{L^2(\bar{C}(9\sigma_1/8))}+1)
    \left[\log(C \|w_\theta(\cdot,t)\|_{L^2(\bar{C}(9\sigma_1/8))}
     +C \|v(\cdot,t)\|_{L^2(\bar{C}(9\sigma_1/8))}+e)\right]^\frac{1}{2}.
\eal
\] Here $a(t)$ is the average of the $r L_\theta(\cdot, t)$  in $\bar{C}(9\sigma_1/8))$
under the 2 dimensional volume element.  After scaling as done at the end of Section 2,
we find, for $k=|x'|$ that,
\[
\al
|L_\theta(x, t)| &\le \frac{C | \ln k |^{1/2}}{ k^{1/2}} +
 C k^{-3} \int_{C_{0.5 k, 5 k}} | L_\theta| dy\\
&\le \frac{C | \ln k |^{1/2}}{ k^{1/2}} +
 C k^{-3} \left(\int_{C_{0.5 k, 5 k}} | L_\theta|^6 dy \right)^{1/6} k^{5/2}\\
&\le \frac{C | \ln k |^{1/2}}{ k^{1/2}} +
 C k^{-1/2} \left(\int_{\reals^3} | L_\theta|^6 dy \right)^{1/6}\\
&\le \frac{C | \ln k |^{1/2}}{ k^{1/2}} +
 C k^{-1/2} \left(\int_{\reals^3} | v |^2 dy \right)^{1/2}.
\eal
\]In the above we just used the 3 dimensional Sobolev inequality and the energy estimate for the
velocity $v$.
Thus $|L_\theta(x, t)| \le \frac{C | \ln k |^{1/2}}{ k^{1/2}}$ for $|x'| = k \le \min \{1/2, R \}$.
\qed

We remark that the energy estimate scales as $-1/2$ power as the distance. Since $L_\theta$
is scaling invariant, the bound in the proposition scales in the same way, modulo the log term.

\section*{Acknowledgements}

Zhen Lei was in part supported by
NSFC (grants No. 11171072, 11421061 and 11222107), SGST
09DZ2272900, Shanghai Talent Development Fund and Shanghai Municipal Dawn project. 
Q. S. Zhang gratefully acknowledges the supports by
Siyuan Foundation through Nanjing University and by Simons
Foundation.
We should also thank Mr. Pan Xinghong for discussions on the
problem. Finally, we are indebted to the anonymous referee who checked the
paper very carefully and found an error in an earlier version of
the paper.

\end{document}